\documentclass[12pt]{article}
\advance\hoffset by -1.1truecm

\def\be{\begin{equation}}
\def\ee{\end{equation}}
\def\ba{\begin{eqnarray}}
\def\ea{\end{eqnarray}}

\newcommand{\x}{{\bf x}}
\newcommand{\p}{{\bf p}}

\newcommand{\sn}{\smallskip\newline}
\newcommand{\mn}{\medskip\newline}

\begin{document}
\title{Three short distance structures from quantum algebras}
\author{Achim Kempf\thanks{Research Fellow of Corpus Christi 
College in the University of Cambridge}\\
Department of Applied Mathematics \& Theoretical Physics\\
University of Cambridge, U.K.}
\date{ }
\maketitle
\begin{abstract}
We review known and we present new results on three types of short 
distance structures of observables which typically appear in studies of
quantum group related algebras. In particular, one of the short distance structures 
is shown to suggest a new mechanism for the introduction of internal symmetries.
\end{abstract} 
\vskip-11truecm

\hskip11.5truecm
{\tt DAMTP/97-100}

\hskip11.5truecm
{\tt q-alg/9709023}
\vskip10.5truecm

\section{Introduction} 
In studies of quantum group related associative algebras, see e.g.
\cite{DRI}-\cite{wessetal}, the 
commutation relation 
\be 
 ab -q ba = 0, \qquad, \qquad q \in {\bf C}
\label{aa1}
\ee
is among the most typical, together with the inhomogenous relation:
\be
\tilde{a} \tilde{b} - q^\prime \tilde{b} \tilde{a} = q^{\prime\prime}, 
\qquad \qquad q^\prime,q^{\prime\prime}\in {\bf C}
\label{aa2}
\ee
As pointed out e.g. in \cite{smbook}, for $q\ne 1$ Eqs.\ref{aa1} and \ref{aa2}
can be transformed into another e.g. through
\be
b = \tilde{b},  \qquad a = \tilde{a} + \tilde{b} q^{\prime\prime}(q-1)^{-1}, \qquad q=q^\prime 
\label{aa3}
\ee
In the following we will focus on Eq.\ref{aa2}.
\sn
For quantum mechanical applications it is crucial to identify 
a canonical pair of hermitean operators, say
$\x$ and $\p$ that play the role of observables with real expectation values. With the 
development of a generalised, noncommutative linear algebra 
(see e.g.
\cite{sm-tianjin} for the braided case) and ultimately functional analysis
it may  prove possible and fruitful to generalise the very concept of hermiticity.
For the time being, however,  we will  have to choose observables which are hermitean
in the conventional sense, i.e. when represented on a separable complex Hilbert space. 
Any choice of observables will
be characterised by a choice of commutation relations. At least on a formal algebraic level,
many choices of commutation relations can be transformed into another, as e.g. the example
of Eqs.\ref{aa1},\ref{aa2},\ref{aa3}
 shows. We also recall that all separable complex Hilbert spaces are isomorphic.
The deeper reason why choices of observables with different commutation relations can, and
indeed must exhibit new phsical features, such as a modified short distance structure, is the fact
that transformations which change commutation relations 
cannot be implemented unitarily, as is not difficult to verify.
\mn
The first ansatz for an identification of observables is to
try to identify $\tilde{a}$ and $\tilde{b}$ directly with observables $\x$ and $\p$, obeying
\be
\x\p-q\p\x =i\hbar
\label{aa4}
\ee
We will consider this approach only briefly.
In this case, if e.g. $\p =\p^\dagger$ is chosen then $\x=\x^\dagger$ is not consistent with Eq.\ref{aa4}
for $q\ne 1$, so that an alternative position observable needs to be defined.  These studies have
been carried out in detail, see e.g. \cite{wessetal}. One of the main results concerning the
resulting short distance structure has been that both, position and momentum space become 
discretised (with exponential spacing).
\mn
Let us now turn to a different identification of observables in Eq.\ref{aa2}, which will yield the second
and third short distance structure which we will here consider. The ansatz is to
take as $\tilde{a}$ and $\tilde{b}$ two mutually adjoint operators $a $ and $a^\dagger$ obeying
(this $a$ is not to be confused with the $a$ of Eq.\ref{aa1})
\be
a a^\dagger - q a^\dagger a = 1, \qquad \qquad q \in {\bf R}
\label{aaa}
\ee
Observables $\x$ and $\p$ are then identified as usual as the 
hermitean and anti-hermitean parts of $a$
\be
\x = L (a+a^\dagger), \qquad \qquad \p = i K (a-a^\dagger)
\ee
where $L$ and $K$ carry units of length and momentum respectively.
The involution defined on the generators as $\x=\x^\dagger, \p=\p^\dagger$ is
consistent, and Eq.\ref{aaa} translates into
\be
[\x,\p] = i\hbar (1 + (q-1)(\x^2/4L^2 + \p^2/4K^2))
\ee
with $KL = \hbar (q+1)/4$. This approach (also for the more general case of 
$U_q(n)$ comodule algebras) has been introduced in \cite{ixtapa,ak-jmp-ucr}, where also
the corresponding short distance structure has been analysed: For $q>1$ both the uncertainty
in position and in momentum are separately finitely bounded from below. The main features
can also be studied in the simpler case of only one correction term on the RHS of the commutation relation:
\be
[\x,\p] = i\hbar(1+\beta \p^2) 
\label{ccr}
\ee
The two cases $\beta >0$ and $\beta <0$ (corresponding to $q >1$ and $q<1$) lead to very different
short distance structures. The case of $\beta >0$ has been studied in considerable detail.
The uncertainty relation then reads
\begin{equation}
\Delta x \Delta p \ge \frac{\hbar}{2} (1+\beta (\Delta p)^2+...) \label{1}
\end{equation}
In the context of string theory and quantum gravity generalised uncertainty relations
and in particular  relations of this type have been discussed, see e.g.
 \cite{townsend}-\cite{camelia}. The uncertainty relation then expresses the presence of a  
natural ultraviolet cutoff as expected in the region of  the Planck scale. Technically,
for all states $\vert \psi \rangle$ in 
a representation of the commutation relations there holds 
$\Delta x_{\vert \psi\rangle} \ge \Delta x_0 = \hbar \sqrt{\beta}$.
We will later return to the case $\beta >0$. 
\mn
Let us now begin with the
analysis of the case $\beta <0$, which is new.
The uncertainty relation in this case
yields for the minimal position uncertainty $\Delta x_0 =0$, as usual.
Taking the trace on both sides of Eq.\ref{ccr} shows
that finite dimensional representations are no longer excluded.
Indeed, there now exist even one-dimensional representations with $\x$
represented as some arbitrary number and $\p$ represented as
$\pm \vert \beta \vert^{-1/2}$. 
All finite dimensional representations reduce to direct sums of 
these: In the $\p$ eigenbasis $\p_{ij}=\p_i\delta_{ij}$
 and the commutation relations, $\x_{rs}(\p_s-\p_r)
=i\hbar \delta_{rs}(1+\beta \p_r^2)$ yield $\p_r =
 \pm \vert \beta\vert^{-1/2}$, thus $([\x,\p])_{rs} = 0$, so that $\x$ is
diagonalisable simultaneously with $\p$, and we obtain $\p_{rs} 
= \mbox{diag}(p_1,p_2,...,p_n)$
and $\x_{rs} = \mbox{diag}(x_1,x_2,...,x_n)$ with $p_i \in 
\{-\vert \beta\vert^{-1/2}, \vert \beta\vert^{-1/2}\}$ and $x_i \in {\bf R}$. 
There also appear new features in 
infinite dimensional representations. Consider
 the spectral representation of $\p$:
\begin{eqnarray}
\p.\psi(\lambda) &=& \lambda \psi(\lambda)\label{e1} \\
\x.\psi(\lambda) &=& i\hbar \left(\frac{d}{d\lambda}+\beta \lambda 
\frac{d}{d\lambda} \lambda\right)
\psi(\lambda)\label{e2}\\
\langle \psi_1\vert \psi_2 \rangle &=& \int_I d\lambda~ \psi_1^*(\lambda)
\psi_2(\lambda)\label{int}
\end{eqnarray}
The family of
 operators $G$ defined through the
 integral kernel ($a,b \in {\bf C}$)
\begin{equation}
G(\lambda,\lambda^\prime) =  \left( a ~\Theta(\lambda -
 \vert \beta\vert^{-1/2}) + b ~\Theta(\lambda + \vert \beta\vert^{-1/2})
\right) \delta(\lambda -\lambda^\prime) 
\end{equation}
commute with both $\x$ and $\p$. Each $G$  is diagonal
 and constant apart from two steps
where it cuts momentum space, and with it the representation,
 into three unitarily inequivalent parts.
The representation which has the proper limit as $\beta
\rightarrow 0$ is given by Eqs.\ref{e1}-\ref{int} 
with the integration interval
$I:= I_c = [-\vert \beta\vert^{-1/2}, \vert \beta\vert^{-1/2}]$. 
Thus, $\p$ becomes a bounded self-adjoint operator. Let us  
calculate the defect indices of $\x$ in this representation, i.e. 
the dimensions of the kernels of $(\x^*\pm i)$. To this end we
check for square integrable solutions to
\begin{equation}
i \left(\partial_\lambda -\vert \beta\vert(\lambda^2 \partial_\lambda
 +\lambda)\right) \psi_\xi(\lambda) = \xi \psi_\xi (\lambda)
\label{pev}
\end{equation}
with $\xi = \pm i$ (from now on we 
set $\hbar =1$). Eq.\ref{pev} is solved by 
\begin{equation}
\psi_\xi(\lambda) = 
\langle \lambda \vert \xi\rangle =
N \left(1-\vert \beta\vert \lambda^2\right)^{-1/2} 
\left(1- \sqrt{\vert \beta\vert} \lambda\right)^{i \xi/2} 
\left(1+\sqrt{\vert \beta\vert} \lambda\right)^{-i \xi/2} 
\label{sols}
\end{equation}
which are non-square integrable on $I_c$ for all $\xi \in {\bf C}$,
 in particular also for $\xi=\pm i$. Thus, 
the defect indices are $(0,0)$, i.e. $\x$ is still essentially 
self-adjoint with a unique spectral representation (compare with
the well known fact that 
the operator $i\partial_\lambda$ which ordinarily represents 
$\x$ on momentum space has defect indices (1,1) on the intervall).
The position eigenfunctions
are given by Eq.\ref{sols} for real $\xi$.
With the continuum normalisation
  $N=(2\pi)^{-1/2}$ it is not difficult to verify 
orthonormalisation and completeness:
\begin{eqnarray}
\int_{-\vert\beta\vert^{-1/2}}^{\vert\beta\vert^{-1/2}}
d\lambda ~ \langle \xi \vert\lambda\rangle
\langle\lambda\vert\xi^\prime
\rangle &=& \delta(\xi-\xi^{\prime})\\
\int_{-\infty}^\infty d\xi~ \langle 
\lambda\vert\xi\rangle\langle\xi\vert
\lambda^\prime\rangle &=& 
\delta(\lambda - \lambda^\prime)
\end{eqnarray}
The generalised Fourier factor given in Eq.\ref{sols}
 yields 
 $\psi(\xi) =  \int_{I_c} d\lambda
\langle \xi\vert\lambda\rangle \psi(\lambda)$ for the  transformation
which maps momentum space wave functions $\psi(\lambda)=
\langle \lambda \vert \psi\rangle$ 
to position space wave functions $\psi(\xi) = \langle 
\xi \vert \psi \rangle$. To summarise in the case of our second short distance structure,
the case of  $\beta <0$, we have found
the presence of a continuous position spectrum while, unexpectedly, 
momentum space becomes bounded.
\medskip\newline
Let us now turn to the case $\beta>0$ and our third short distance structure.
As is well known, and as is easily derived from Eq.\ref{1},
the position resolution $\Delta x$ now becomes finitely
 bounded from below: $\Delta x_0= \sqrt{\beta}$. 
To be precise,
for all normalised vectors $\vert \psi\rangle $
in a domain $D$ on which the 
commutation relations hold the
position uncertainty obeys $\Delta x_{\vert\psi\rangle}
 =\langle \psi \vert (\x -\langle\psi\vert \x
\vert\psi\rangle)^2\vert\psi\rangle^{1/2} \ge 
\Delta x_0 =\sqrt{\beta} $. A convenient representation is given by 
Eqs.\ref{e1}-\ref{int} with $I= {\bf R}$.
On any dense domain $D$ in a Hilbert 
space $H$
 on which the commutation relations hold the position 
operator can only be symmetric
 but not self-adjoint, as diagonalisability is
 excluded by the
 uncertainty relation (eigenvectors to an observable  
automatically have
vanishing uncertainty in this observable). This
 also excludes the possibility of finite 
dimensional representations of the
 commutation relations (since in these symmetry
 and self-adjointness coincide), 
as could of course also be seen by taking the 
trace of both sides of Eq.\ref{ccr}. 
The underlying
functional analytic structure was first discussed in
\cite{ixtapa,ak-jmp-ucr}. 
\medskip\newline
We will now discuss further physical 
implications of this short distance structure, related to internal symmetries. 
Our aim is to show that 
the unobservability of localisation beyond 
the minimal uncertainty $\Delta x_0$ 
represents
an internal symmetry where degrees of freedom which
 correspond to small scale structure beyond the
Planck scale turn into internal degrees of freedom.

Consider a 
$*$-representation (such as given by Eqs.\ref{e1}-\ref{int})
 of the commutation relation
 Eq.$\ref{ccr}$ with $\beta >0$
on a maximal dense domain $D$ in a Hilbert space $H$. 
Then $\x$ is merely symmetric, i.e. $D$ is smaller
 than the domain $D_{\x^*}$ of the 
adjoint operator $\x^*$ (which is not symmetric). 
The deficiency 
spaces $L_+,L_-$, i.e. the spaces
spanned by eigenvectors of $\x^*$ with eigenvalues
 $+i$ and $-i$ are one-dimensional, i.e. the 
deficiency indices are (1,1). 

Thus, there exists a set of 
self-adjoint extensions of $\x$ which is in
one-to-one correspondence with the set of unitary
 transformations $\tilde{U} : L_+ \rightarrow L_-$. 
We recall that, by the usual procedure, each $\tilde{U}$
defines a unitary extension of the Cayley transform of
 $\x$, with the inverse Cayley transform then defining 
a self adjoint extension of $\x$. On the eigenvalues,
 Cayley transforms are M{\"o}bius transforms.

 The $\tilde{U}$ differ exactly by
the set $G$ of unitary transformations $U:
 L^+ \rightarrow L^+$, which we may here call the `local group' $G$.
Thus, the set of self-adjoint extensions $\{\x_\alpha\}$ 
forms a representation of the local group, where $\alpha$
 which labels the
self-adjoint extensions is a
vector in the fundamental representation of $G$.
The local group also acts on the
 set of spectra $\{\sigma_\alpha\}$
of the $\x_\alpha$. Let us denote the eigenvalues of the self-adjoint
extension $\x_\alpha$ by $v_\alpha(r)$. Then, for any fixed $r$,
we obtain an orbit $O(r):= \{ v_{U.\alpha}\vert U\in G\}$
 of eigenvalues under the action of $G$. 

For the case of Eq.\ref{ccr} 
the scalar product of eigenvectors of $\x^*$ has been
calculated in \cite{kmm}:
\begin{eqnarray}
\langle \xi \vert \xi^\prime \rangle  = 
 \frac{2 \sqrt{\beta}}{
\pi (\xi - \xi^{\prime})} {\mbox{ }} \sin\left(\frac{\xi 
-\xi^{\prime}}{2 \sqrt{\beta}} 
\pi\right)
\label{spfp}
\end{eqnarray}
{}From its zeros
 we can read off the family of discrete
 spectra of the self-adjoint extensions:
\begin{equation}
\sigma_\alpha = \left\{v_{\alpha}(r) = 
(2 r + s/\pi) \sqrt{\beta}~\vert~
r \in \mbox{N}\right\} ~~~~~\mbox{where}~~~~~
\alpha = e^{is} \mbox{~~ with ~~}s \in [0,2\pi[
\end{equation} 
The spectra are equidistant and self-adjoint
 extensions
 differ by a shift of their lattice of eigenvalues.
 The local group is here
the group of translations of the lattices of eigenvalues. 
Due to the periodicity of the lattice this group is 
$S^1$, or $U(1)$. This had to be expected since in this case
$L_+$ is one-dimensional and the
 self-adjoint extensions therefore form a representation 
of the local group U(1). 
 
Each choice of self-adjoint extension of the 
position operators therefore
corresponds to a choice of lattice on which 
the physics takes place. However, the commutation
 relations also
imply that the smallest uncertainty in positions 
becomes finite and large enough so that the actual choice
of lattice cannot be resolved. Technically, all
self adjoint extensions of $\x$ coincide when 
restricted to a domain $D$ on which the commutation relations hold.

If, therefore,
 with a physical state $\vert \psi\rangle \in D$ 
also some vector $\alpha$ is specified, as a
choice of self-adjoint extension, the action should 
be invariant. We therefore arrive at a global symmetry principle where
the additional information given by the `isospinor' $\alpha$ can be 
interpreted. Assume that the state of a particle is
 projected onto 
a state of maximal localisation ($\Delta x = \Delta x_0$)
with position expectation $\xi$. Specifying $\alpha$ is
 to specify one point in the orbit of
the eigenvalue $\xi$ under the action of the local group.
 As a convention one can specify that this is  
where the maximally localised particle (or one point of
 it if it is viewed as
an extended particle) is said to "actually" sit.
This is consistent because    
the radius of the orbits of the eigenvalues 
is $\sqrt{\beta} = \Delta x_0$ i.e. of the size of
 the finite minimal uncertainty $\Delta x_0$, so that
 all these conventions,
differing only by the action of the local group, 
cannot be  distinguished observationally. 
For example, the pointwise multiplication of fields
 as discussed
e.g. in \cite{kmm} can be reformulated in terms of
 a choice of position eigenbasis, 
rather than the set of maximally localised fields. 
The gauge principle is that the action is invariant under
the local group. We remark that the proof of ultraviolet 
regularity still goes through since not only the
fields of maximal localisation, but also the
 position eigenfields are normalisable. The `local' group may then
also be taken to act
locally, i.e. we consider $\vert \psi\rangle \in D \otimes L_+$.
It is unobservable whether one specifies one
 self-adjoint extension's lattice here and another's
there, as long as the parallel transport of
 $\alpha$ is consistently defined. At large scales this should 
turn into the ordinary local gauge principle.

There is therefore a possibility that internal
symmetry spaces arise as 
deficiency spaces of position operators.
 Introducing $\Delta x_0>0$ 
the infinite dimensional Hilbert space of
 fields develops
 special dimensions that
correspond to degrees of freedom that describe localisation
beyond what can be resolved, and which can therefore be 
viewed as internal degrees of freedom. 
In the case of one dimension our studies yield the
simple intuitive picture that certain 
corrections to the 
uncertainty relations lead to physics 
on a whole set of possible
lattices, while the choice of any 
particular lattice from the set
cannot be resolved and does therefore
 correspond
 to an internal degree of freedom. 
This may be a new mechanism, or it could be 
related to the Kaluza Klein
mechanism. It should be very interesting to explore possible connections
to string theory. 
\mn
In particular, it should be interesting to explore a possible connection to an
observation by Faddeev: In \cite{Faddeev} he noted a simple mechanism by which
in the process of compactification or discretisation one degree of freedom with
infinitesimal generators
of $\x$ and $\p$ with $[\x,\p]=i\hbar$ can turn into two degrees of freedom
of $\xi_1,\pi_1$ and $\xi_2, \pi_2$ where the $\xi_i$ and $\pi_i$ are now finite
translations and boosts:  $\xi_i=\exp(ia_i\x), ~\pi_i=\exp(b_i\p)$, with appropriately
chosen parameters $a_i$ and $b_i$, yielding a modular group
action. These degrees of freedom would here correspond to the (external) degree of freedom
which measures distances in units of the fundamental length, and the (internal) degree
of freedom which is related to the position of the lattice of eigenvalues, 
technically analogous to the
pair of winding and ordinary modes degrees of freedom. 
\mn
Finally, let us remark that general commutation relations, for all types of
short distance structure, as long as $\x_i=\x_i^\dagger$ and $\p_i=\p_i^\dagger$, 
can  be introduced into the quantum field
 theoretical path integral. The main observation is that 
the functional analysis of representations of the commutation
 relations on wave
 functions extends to representations on fields. 
For early considerations on the representation of the $\x,\p$
commutation relations on the space of fields which is formally
being summed over in the field theoretic path integral see e.g.
\cite{deWitt}. In particular, 
$\Delta x_0 = \hbar \sqrt{\beta} >0$, in spite of being
an ensemble-cutoff rather than an individual-case-cutoff,
has been shown to regularise the ultraviolet 
in euclidean field theory. The issue of regularisation 
through $\Delta x_0>0$ has been
studied extensively in \cite{ft}-\cite{km}.
\medskip\newline
Our analysis concerning internal symmetries
covered  the case of one 
space-like dimension with a positive correction term to 
the commutation relations. 
Due to the indefinite Minkowski signature the
temporal coordinate can be  expected to come with correction terms
with a negative sign i.e. with a short distance structure of the second
type ($\beta <0$) which we here considered. Our study above of this case 
indicates that this direction,
i.e. $\x_0$, should have a unique self-adjoint extension and that it therefore 
may not contribute to the internal symmetries. 
The generic case of  $d$ dimensions with Minkowski signature
and possibly higher order correction terms will require careful 
investigation. It is of course possible that more than the here covered three short
distance structures may appear.

\end{document}